\documentclass[a4paper,10pt]{article}

\usepackage[english]{babel}
\usepackage[latin9]{inputenc}
\usepackage{amsmath}
\usepackage{amsfonts}
\usepackage{amssymb}
\usepackage{amsthm}
\usepackage{graphics}
\usepackage{fullpage}
\usepackage{enumerate} 
\usepackage{multicol} 
\usepackage{xy} 
\xyoption{all}
\usepackage{lscape}
\usepackage{nomencl}
\usepackage{pdfpages}
\usepackage{hyperref}

\hypersetup{
    backref=true,    
    pagebackref=true,
    hyperindex=true, 
    colorlinks=true, 
    breaklinks=true, 
    urlcolor= red, 
     linkcolor= black, 
     citecolor=black,
     bookmarks=true, 
     bookmarksopen=true,            
}

\newcommand{\fct}[5]{\raisebox{-.25cm}{ 
\begin{tabular}{rccc}
    #1 : & #2 & $\rightarrow$ & #3 \\
   & #4 & $\mapsto$ & #5 \\
\end{tabular}}}
\newcommand{\ie}{, \textit{i.e.} }
\newcommand{\R}{\mathbb{R}}
\newcommand{\C}{\mathbb{C}}
\newcommand{\Z}{\mathbb{Z}}
\newcommand{\N}{\mathbb{N}}
\newcommand{\K}{\mathbb{K}}
\newcommand{\g}{\mathfrak{g}}
\newcommand{\bb}[1]{\mathbb{#1}}
\newcommand{\ca}[1]{\mathcal{#1}}
\newcommand{\mk}[1]{\mathfrak{#1}}
\newcommand{\id}{\mathop{\rm Id}}
\newcommand{\der}[1]{\mathop{\rm Der}(#1)}

\newcommand{\ad}[1]{\mathop{\rm ad}\left(#1\right)}

\newcommand{\pt}{\le}
\newcommand{\aut}[1]{\mathop{\rm Aut}(#1)}
\newcommand{\ssi}{if and only if\ }
\newcommand{\exi}{, there exists\ }
\newcommand{\prtt}{for all\ }
\newcommand{\tr}{\mathop{\rm Tr}}
\newcommand{\Hom}{\mathop{\rm Hom}}
\renewcommand{\hom}[1]{\mathop{\rm Hom}^{(i)}\left(#1\right)}
\newcommand{\End}[1]{\mathop{\rm End}(#1)}

\newcommand{\Gl}{\mathop{\rm Gl}}

\newcommand{\fibre}[2]{\raisebox{.55cm}{$\xymatrix{#1 \ar@{->>}[d] \\ #2}$}}

\renewcommand{\u}{\mk{u}\!}
\newcommand{\U}{U\!}

\theoremstyle{definition}
\newtheorem{defi}{\underline{Definition}}[section]
\theoremstyle{plain}
\newtheorem{theo}{\underline{Theorem}}[section]
\newtheorem{prop}{\underline{Proposition}}[section]
\newtheorem{lem}{Lemma}[section]
\newtheorem{coro}{\underline{Corollary}}[section]

\theoremstyle{remark}
\newenvironment{rem}{\underline{\emph{Remark}} :\ }{}
\newenvironment{demo}{\begin{proof}

}{\end{proof}}

\title{Generalized flag geometries associated with $(2k+1)$-graded Lie algebras\\[.5cm]}
\author{Chenal, Julien\\ \tiny{Institut Elie Cartan Nancy(IECN),}\\[-.2cm] \tiny{Nancy-Universit\'e, CNRS, INRIA,}\\[-.2cm] \tiny{Boulevard des Aiguilettes B.P. 239}\\[-.2cm] \tiny{F-54506 Vandoeuvre-l\`es-Nancy}\\[-.2cm]\tiny{chenal@iecn.u-nancy.fr}}
\date{}
\begin{document}


\maketitle

\medskip

\medskip

\begin{abstract}

In this paper, we present the construction of a geometric object, called a \emph{generalized flag geometry} $(X^+, X^-)$, corresponding to a $(2k+1)$-graded Lie algebra $\g=\g_k\oplus\dots\oplus\g_{-k}$. We prove that $(X^+, X^-)$ can be realized inside the space of \emph{inner filtrations} of $\g$ and we use this realization to construct ``algebraic bundles'' on $X^+$ and $X^-$ and some sections of these bundles. Thanks to these constructions, we can give a realization of $\g$ as a Lie algebra of polynomial maps on the positive part of $\g$, $\mk{n}^+_1:=\g_1\oplus\dots\oplus\g_k$, and a trivialization in $\mk{n}^+_1$ of the action of the group of automorphisms of $\g$ by ``birational''maps.\\[.3cm]

\end{abstract}

\section{Introduction}

\

It is well-known (see \cite{Loos75}) that \emph{$3$-graded Lie algebras}
$\g = \g_1 \oplus \g_0 \oplus \g_{-1}$ give rise to \emph{Jordan pairs}
$(V^+,V^-)=(\g_1,\g_{-1})$, where the trilinear Jordan-structure is given by
$T^\pm(x,y,z)=[[x,y],z]$, and that, conversely, the $3$-graded Lie algebra
can be reconstructed from a (linear) Jordan pair by the \emph{Tits-Kantor-Koecher construction}.
The Jordan pair can then be interpreted as the ``infinitesimal structure'' corresponding
to a global geometry $(X^+,X^-)$, called \emph{generalized projective geometry} in
\cite{B02}.
Thus $3$-graded Lie algebras correspond to generalized projective geometries,
 and in \cite{BN1}, an explicit and purely algebraic construction of a ``universal model''
of these geometries in terms of \emph{$3$-filtrations} of the Lie algebra has been given.

\

The topic of the present work is to extend these results to $\Z$-graded Lie algebras of the type 
$$\g = \g_k \oplus \g_{k-1} \oplus \ldots \oplus \g_{-k},$$
which we call \emph{$(2k+1)$-graded}. Especially the $5$-graded case has already attracted much attention, mainly because
every simple complex Lie algebra (of dimension bigger than $3$) admits a $5$-grading
having the special feature that the spaces $\g_{\pm 2}$ are one-dimensional, and then
the pair $\g_{\pm  1}$ corresponds to \emph{Freudenthal-Kantor pairs} (see \cite{Kan72}, \cite{Kam} or \cite{AllisFaul}); also, this particular situation is related to geometric topics such as
\emph{symplectic geometry} (see \cite{Rub}) and
\emph{quaternionic symmetric spaces} (the so-called Wolf-spaces, cf. \cite{B02q}).
The typical situation that we consider here is in a certain sense opposite to this very particular
case: we do not assume $\g$ to be simple, nor finite-dimensional, nor even defined over
a field but rather over a commutative base ring $\K$ (the only assumption we will need is
that certain integers are invertible in $\K$). In this general situation, it is much less clear
what object replaces the Jordan pair from the $3$-graded case; it is known, by work of
I. Kantor, that the pair $(\g_1,\g_{-1})$ is a \emph{generalized Jordan pair} (a sort of
non-commutative version of a Jordan pair, see \cite{Kan70} and \cite{Kan72}), but this pair seems to be
``too small'' if one is looking for an object corresponding to the geometry of
the whole graded algebra. Indeed, the problem of finding the good analog of the Jordan
pair in the higher graded case remains open, but we hope that a good understanding of
the geometry $(X^+,X^-)$ we construct  will be a step towards the solution of this problem. 

\

The geometry associated to a $(2k+1)$-graded Lie algebra $\g$, which we call a \textit{generalized flag geometry}, is constructed as follows:
we consider the space $\cal F$ of all \emph{$(2k+1)$-filtrations of $\g$}, and we define
a natural notion of \emph{transversality of filtrations} (Definition \ref{deffil}).
We prove that two filtrations $e$ and $f$
are transversal (denoted by $e \top f$) if and only if they come from a grading
(then $e$ is the ``ascending'' and $f$ the ``descending'' filtration, see Part \ref{transverseetgraduation}) of Theorem \ref{theoprincipal}).
Our key result (Part \ref{Uagitsimplementtransitivement}) of Theorem \ref{theoprincipal}) says that the set $f^\top$ of all filtrations
$e$ that are transversal to a given filtration $f$ carries a natural \emph{torsor structure}
(where we use the term ``torsor'' in the same way as in \cite{BK1}):
there is a natural, simply transitive group action on $f^\top$, hence, after choice of
an arbitrar origin, $f^\top$ becomes a group. In the case of $3$-filtrations, this group is
a vector group (underlying an affine space, see \cite{BN1}), whereas for $k>1$, it is no longer
abelian but rather unipotent, with ``Lie algebra'' the positive (or negative, according to choice) part
$\mk{n}^+_1:=\g_1 \oplus \ldots \oplus \g_k$, and instead of an affine space, we have a
dilatation action making this group a \emph{Carnot group} (see, e.g., \cite{Bul07}).
These facts make the higher graded situation geometrically much more complicated, as we
shall explain now.

\

We fix a $(2k+1)$-grading of $\g$ as ``base point'', which corresponds to two transversal
filtrations $o^+$ and $o^-$, and we consider the ``connected components''
$X^\pm$ of $\cal F$ containing $o^\pm$ (Section \ref{varietededrapeaux}). For any $f \in X^\pm$, the torsor
$f^\top$ belongs to $X^\mp$ (Theorem \ref{theostructure}); we call this structure a \emph{linear pair geometry}.
Our aim is then to describe geometric structures of  $X^\pm$ with respect to
the ``charts'' $f^\top$: first of all, we give an algebraic description of the ``intersection of
chart domains'' $f^\top \cap e^\top$ in terms of certain cocycles called \emph{denominators}
and \emph{co-denominators} (Definition \ref{defdeno} and Corollary \ref{corog.xinV+}).
These are related to a generalization of the well-known \emph{Bergman operator} from
Jordan-theory; more intrinsically, we define, in a purely algebraic way, certain
vector bundles over $X^\pm$, and then relate these operators to invariant sections of such
bundles,  called the \emph{canonical kernel} (Definition \ref{noyaucano} and Theorem \ref{noyaucaniosection}).
Next, we describe the action of the \emph{projective elementary group} (which is the natural
symmetry group in this context, acting transitively on $X^\pm$, see Section \ref{groupprojelem}) with respect to
the ``charts'' $f^\top$. This description is explicit and algebraic (in case of finite dimension
over a field, it is indeed a birational group action), but the formulas are considerably more
complicated than in the $3$-graded case (Proposition \ref{actionrationnelle}).
They generalize the classical description of the projective group on projective space,
given in an affine picture by  fractional-linear maps (see Example \ref{actionparhomographie}).
In subsequent work, we will use this description for proving that, under suitable 
topological assumptions on $\K$ and on $\g$, the geometries $X^\pm$ are indeed
smooth manifolds and that the algebraically defined bundles are smooth bundles with
smooth section given by the canonical kernel. These results will give new and interesting
examples of smooth infinite dimensional geometries, associated to infinite dimensional Lie algebras.

\

As an intermediate step to the description of the group action, we give a realization of
its ``Lie algebra'' (which is essentially $\g$) with respect to the ``chart'' $f^\top$. It
 is given by vector fields (sections of the algebraically defined tangent bundle)
which turn out to be \emph{polynomial} (Section \ref{champspoly}),
but again this description is much more complicated than in the $3$-graded
case; in this case, the polynomials are quadratic, and the algebraic structure of
this polynomial Lie algebra is directly related to the Jordan pair structure (see \cite{BN1} and \cite{Up}). In the higher graded case, the missing analog of the
Jordan pair makes it much more difficult to understand the polynomials arising
from the polynomial realization of $\g$. It seems likely that a
conceptual understanding of this situation should use the notion of \emph{generalized contact
structure}, proposed in \cite{CKR}: indeed, in loc. cit. it is
shown that (under suitable assumptions on a simple Lie algebra over $\R$), the projective elementary group can be
indeed characterized as the group preserving the generalized contact structure, which
thus appears to be a key feature of this situation.
We intend to come back to this question in subsequent work.

\

This results are part of the author's thesis \cite{moi2} and some of them have been anounced in a note \cite{moi1}.

\vspace{.5cm}

\underline{Notation}: In this paper, $\K$ is a commutative ring with $\N\subset \K^{\times}$. Let $k\in\N$.

\section{Filtered and graded Lie algebras}\label{filtrationetgraduation}

\

In this part, we define the two notions of \emph{filtrations} and \emph{gradings} for a $\K$-module and for a Lie algebra and focus on links between them. Let $V$ be a $\K$-module and $\g$ a Lie algebra over $\K$.

\subsection{Filtrations and gradings of $\K$-modules}

\begin{defi}
A \emph{$k$-filtration} of the module $V$ is a flag of subspaces of $V$:$$\mk{f}=\left(0\subset\mk{f}_k\subset\mk{f}_{k-1}\subset\dots\subset\mk{f}_1=\mk{g}\right).$$
Two filtrations $\mk{e}$ and $\mk{f}$ are called \emph{transversal}, and we write $\mk{e}\top\mk{f}$, if $$V=\mk{e}_n\oplus\mk{f}_{k-n+2}, \text{ for } 2\pt n\pt k-1.$$
\end{defi}

\begin{defi}
A \emph{$k$-grading} of $V$ is a family $(V_1,\dots,V_k)$ of subspaces of $V$, such that $$V=V_1\oplus\dots\oplus V_k.$$
\end{defi}
Each grading $(V_1,\dots,V_n)$ defines two transversal filtrations $\mk{f}^+$ and $\mk{f}^-$, defined, for $1\pt n\pt k$, by:$$\begin{array}{rcl} \mk{f}^+_n & := & V_k\oplus\dots\oplus V_n,\\ \mk{f}^-_{k-n+1} & := & V_{n+1}\oplus\dots\oplus V_1.\\ \end{array}$$
Conversely, by \cite{BL08}, two transversal filtrations come from a grading of $V$:
\begin{prop}\label{propBL}
Two filtrations $\mk{e}$ and $\mk{f}$ of $V$ are transversal if and only if they come from a grading of $V$\ie if and only if there exists a grading $V=V_1\oplus\dots\oplus V_k$ such that $\mk{e}=\mk{f}^+$ and $\mk{f}=\mk{f}^-$. Moreover, if $\mk{e}$ and $\mk{f}$ are transversal, the corresponding grading is defined by $V_n:=\mk{e}_n\cap\mk{f}_{k-n+1}$.
\end{prop}
\begin{demo}
See Proposition 3.3 of \cite{BL08}.
\end{demo}

\subsection{Filtrations of Lie algebras}

\

\begin{defi}\label{deffil}
A \emph{$(2k+1)$-filtration} of a Lie algebra $\g$ is a flag of subspaces of $\g$: $$\mk{n}=(0=\mk{n}_{k+1}\subset\mk{n}_k\subset\dots\subset\mk{n}_{-k+1}\subset\mk{n}_{-k}=\g), \text{ such that } \left[\mk{n}_m,\mk{n}_n\right]\subset\mk{n}_{m+n}.$$
Such a flag will be denoted by $\mk{n}=(\mk{n}_k\subset\dots\subset\mk{n}_{-k+1})$.\\
Two $(2k+1)$-filtrations $\mk{m}$ and $\mk{n}$ are called \emph{transversal}, and we write $\mk{m}\top\mk{n}$, if they are transversal flags in the sense of the preceding subsection\ie $\g=\mk{m}_n\oplus\mk{n}_{-n+1}$, for all $-k+1\pt n\pt k$.\\ We say that a filtration is \emph{complemented} if it admits a transversal filtration, and the set of such complemented filtrations will be denoted by $\widetilde{\ca{F}}$. 
\end{defi}
Given a filtration $\mk{n}$, we define the Lie algebra $$\u\left(\mk{n}\right):=\left\{X\in\der{\g}|\ X\mk{n}_n\subset\mk{n}_{n+1},  \text{ for all } -k+1\pt n\pt k\right\}.$$
The elements $X$ of $\u\left(\mk{n}\right)$ are nilpotent. Hence, since integers are assumed to be invertible in $\K$, we can consider the unipotent ``Lie group'' corresponding to $\u\left(\mk{n}\right):$
$$\U\left(\mk{n}\right):=e^{\mk{u}\left(\mk{n}\right)}:=\left\{e^X| X\in\u\left(\mk{n}\right)\right\}, \text{ where } e^X=\sum_{n=0}^{2k}\frac{X^n}{n!}\in\aut{\g}.$$

\subsection{Short $\Z$-gradings of Lie algebras}

\begin{defi}
A \emph{$(2k+1)$-grading} of $\g$ is a family $\left(\g_{-k},\dots,\g_k\right)$ of submodules of $\g$ such that $$\g=\bigoplus_{n=-k}^k\g_n, \text{ with } \left[\g_m,\g_n\right]\subset\g_{m+n} \text{ and } \g_n=0 \text{ if } |n|>k.$$
\end{defi}
If $\g=\g_k\oplus\dots\oplus\g_{-k}$ is $(2k+1)$-graded, the map $D:\g\rightarrow \g$ defined, for $x\in\g_n$, by $Dx=nx$, is a derivation, called the \emph{characteristic derivation} of the grading, which satisfies \begin{equation} \left(D-k\id\right)\dots\left(D-\id\right)D\left(D+\id\right)\dots\left(D+k\id\right)=0.\label{equation}\end{equation} Conversely, any derivation $D$ satisfying (\ref{equation}) is diagonalizable and the eigenspace decomposition is a $(2k+1)$-grading of $\g$. Hence, we can identify the space of $(2k+1)$-gradings of $\g$ with the set $$\widetilde{\ca{G}}:=\left\{D\in\der{\g}|  \left(D-k\id\right)\dots\left(D-\id\right)D\left(D+\id\right)\dots\left(D+k\id\right)=0\right\}.$$
If an element $D\in\widetilde{\ca{G}}$ can be written $D=\ad{E}$, $E\in\g$, then the grading is called \emph{inner}, and the element $E$ is called an \emph{Euler operator} of the grading. We denote by $$\ca{G}:=\left\{\ad{E}|  E\in\g, \ad{E}\in\widetilde{\ca{G}}\right\},$$ the set of inner $(2k+1)$-gradings of $\g$.

\subsection{Transversality and gradings}\label{Transversalityandgradings}

\

If $\g=\g_k\oplus\dots\oplus\g_{-k}$ is a $(2k+1)$-graded Lie algebra, two filtrations are defined by $$\mk{n}^+(D):=(\g_k\subset\g_k\oplus\g_{k-1}\subset\dots\subset\g_k\oplus\dots\oplus\g_{-k+1}),$$
$$\text{and } \mk{n}^-(D):=(\g_{-k}\subset\g_{-k}\oplus\g_{-k+1}\subset\dots\subset\g_{-k}\oplus\dots\oplus\g_{k-1}).$$ The filtration $\mk{n}^+(D)$ is called the \emph{plus-filtration} and $\mk{n}^-(D)$ the \emph{minus-filtration} of the grading. If $D=\ad{E}$ is inner, the filtration $\mk{n}^+(\ad{E})$ is also called \emph{inner} and we denote by $\ca{F}$ the set of all inner filtrations. Observe that, by construction, the filtrations $\mk{n}^+(D)$ and $\mk{n}^-(D)$ defined by a grading of $\g$ are transversal. Conversely, the following holds:

\begin{theo}\label{theoprincipal}
\begin{enumerate}[(1)]
\item\label{transverseetgraduation} Let $\mk{m}$ and $\mk{n}$ be two filtrations of $\g$. Then $\mk{m}$ and $\mk{n}$ are tranversal if and only if they come from a grading of $\g$, defined by $\g_n=\mk{m}_n\cap\mk{n}_{-n}$\ie there exists $D\in\widetilde{\ca{G}}$ such that $\mk{m}= \mk{n}^+(D)$ and $\mk{n}=\mk{n}^-(D)$. In other words, $\widetilde{\ca{F}}=\left\{\mk{n}^+(D), D\in\widetilde{\ca{G}}\right\}$.
\item\label{orbitedeDsousU} Let $\mk{n}=\mk{n}^+(D)\in\widetilde{\ca{F}}$. Then $\U\left(\mk{n}\right)\cdot D=D+\u\left(\mk{n}\right)$, where $\U\left(\mk{n}\right)$ acts  on $\der{\g}$ by $g\cdot D:=g\circ D\circ g^{-1}$.
\item\label{Uagitsimplementtransitivement} Let $\mk{n}\in\widetilde{\ca{F}}$. Then the group $\U\left(\mk{n}\right)$ acts simply transitively on $\mk{n}^\top:=\left\{\mk{m}\in\widetilde{\ca{F}}, \mk{m}\top\mk{n}\right\}$.
\end{enumerate}

\end{theo}

\begin{demo}
\begin{enumerate}[(1)]
\item We have already remarked that $\mk{n}^+(D)\top\mk{n}^-(D)$. Conversely, let us assume that $\mk{m}\top\mk{n}$. According to Proposition \ref{propBL}, the flags $\mk{m}$ and $\mk{n}$ come from a grading of the module $\g$, defined by $\g_n=\mk{m}_n\cap\mk{n}_{-n}$. This is actually a grading in the sense of Lie algebras\ie $\left[\g_m,\g_n\right]\subset\g_{m+n}$. Indeed, $$\left[\g_m,\g_n\right]\subset\left[\mk{m}_m,\mk{m}_n\right]\cap\left[\mk{n}_{-m},\mk{n}_{-n}\right]\subset\mk{m}_{m+n}\cap\mk{n}_{-(m+n)}=\g_{m+n}.$$ Hence, $\mk{m}$ and $\mk{n}$ come from a $(2k+1)$-grading of $\g$.
\item Since $e^X=\text{Id}+\sum\limits_{n=1}^{2k}X^n$, for $X\in\u\left(\mk{n}\right)$, the inclusion $\U\left(\mk{n}\right)\cdot D\subset D+\u\left(\mk{n}\right)$ holds.\\
In order to prove the converse, let us fix $X\in\u\left(\mk{n}\right)$ and show, by induction, that for $2\pt n\pt 2k+1$, there exist $Y_n,R_n\in\u\left(\mk{n}\right)$ such that $$e^{Y_n}De^{-Y_n}=D+X+R_n \text{ and } R_n\left(\mk{n}_i\right)\subset \mk{n}_{i+n}, \text{ for all } -k+1\pt i\pt k.$$
This implies that, for $n=2k+1$, $R_n\left(\g\right)\subset\mk{n}_{k+1}=0$, hence $e^{Y_{2k+1}}De^{-Y_{2k+1}}=D+X$, and the claim follows: $$\U\left(\mk{n}\right)\cdot D\supset D+\u\left(\mk{n}\right).$$
Now, we start the inductive proof: let $n=2$ and define $Y_2:=-X\in\u\left(\mk{n}\right)$. Then $$\begin{array}{rcl} e^{Y_2}De^{-Y_2} & = & \displaystyle{\left(1-X+\sum_{j=2}^{2k+1}(-1)^j\frac{X^j}{j!}\right)D\left(1+X+\sum_{j=2}^{2k+1}\frac{X^j}{j!}\right)}\\
& = & D+\left[D,X\right]+\widetilde{R_2}, \text{ with } \widetilde{R_2}\left(\mk{n}_i\right)\subset\mk{n}_{i+2}.\\
\end{array}$$ Let us consider $\g=\bigoplus\g_n$, the $(2k+1)$-grading of $\g$ defined by $D$. On the one hand, for $x\in\g_i$, $XDx=iXx\in\mk{n}_{i+1}$, while, on the other hand, $DXx=(i+1)Xx+Z_{i+2}$, where $Z_{i+2}\in\mk{n}_{i+2}$. Thus, $\left[D,X\right]x=Xx+Z_{i+2}$, and $$e^{Y_2}De^{-Y_2}=D+X+R_2, \text{ with } R_2\left(\mk{n}_i\right)\subset\mk{n}_{i+2}.$$ Moreover, $R_2=e^{Y_2}De^{-Y_2}-D-X\in\u\left(\mk{n}\right)$, and our claim is proved for $n=2$.\\
Now, let $n\in\bb{N}$ and assume that there exist $Y_n,R_n\in\u\left(\mk{n}\right)$ such that $$e^{Y_n}\cdot D=D+X+R_n \text{ and } R_n\left(\mk{n}_i\right)\subset\mk{n}_{i+n}.$$
We define $Y_{n+1}$ by $$Y_{n+1}:=Y_n+\frac{1}{n}R_n\in\u\left(\mk{n}\right).$$
It follows that $$\begin{array}{rcl} e^{Y_{n+1}} & = & \displaystyle{e^{Y_n+\frac{1}{n}R_n}=\sum_{j=0}^{2k+1}\frac{\left(Y_n+\frac{1}{n}R_n\right)^j}{j!}}\\
& = & \displaystyle{1+Y_n+\frac{1}{n}R_n+\sum_{j=2}^{2k+1}\frac{\left(Y_n+\frac{1}{n}R_n\right)^j}{j!}}\\
& = & \displaystyle{e^{Y_n}+\frac{1}{n}R_n+\widetilde{R^+_{n+1}}}, \text{ with } \widetilde{R^+_{n+1}}\left(\mk{n}_i\right)\subset\mk{n}_{i+n+1}.\\ \end{array}$$
Similarly, $e^{-Y_{n+1}}=e^{-Y_n}-\frac{1}{n}R_n+\widetilde{R^-_{n+1}}$, with $\widetilde{R^-_{n+1}}\left(\mk{n}_i\right)\subset\mk{n}_{i+n+1}$.\\
For $x\in\g_i$, $De^{-Y_{n+1}}x=De^{-Y_n}x-\dfrac{i+n}{n}R_nx+Z^{(1)}_{i+n+1}, \text{ with } Z^{(1)}_{i+n+1}\in\mk{n}_{i+n+1}$. Hence
$$\begin{array}{rcl} e^{Y_{n+1}}De^{Y_{n+1}}x & = & e^{Y_n}De^{-Y_n}x-\dfrac{i+n}{n}e^{Y_n}R_nx+\dfrac{1}{n}R_nDe^{-Y_n}x+Z^{(2)}_{i+n+1}, \text{ with } Z^{(2)}_{i+n+1}\in\mk{n}_{i+n+1}\\
& = & Dx+Xx+R_nx-\dfrac{i+n}{n}R_nx+\dfrac{i}{n}R_nx+Z^{(3)}_{i+n+1}, \text{ with } Z^{(3)}_{i+n+1}\in\mk{n}_{i+n+1}\\
& = & Dx+Xx+Z^{(3)}_{i+n+1},\\
\end{array}$$
thus $e^{Y_{n+1}}\cdot D=D+X+R_{n+1}$, with $R_{n+1}\left(\mk{n}_i\right)\subset\mk{n}_{i+n+1}$. As explained above, the claim follows.
\item First, we prove the transitivity of the action. Let $\mk{n}=\mk{n}^+(D)$ be a filtration and $\mk{m}\in\mk{n}^{\top}$ a filtration transversal to $\mk{n}$. According to Part \ref{transverseetgraduation}), there exists a grading $D'$ such that $\mk{m}=\mk{n}^-(D')$ and $\mk{n}=\mk{n}^+(D')$. We claim that $D-D'\in\u\left(\mk{n}\right)$. Indeed, if we denote by $$\g=\bigoplus_{n=-k}^k\g_n\text{ and } \g=\bigoplus_{n=-k}^k\g'_n$$ the gradings defined by $D$ and $D'$, as $\mk{n}^+(D)=\mk{n}=\mk{n}^+(D')$, we have $\g_k\oplus\dots\oplus\g_n=\mk{n}_n=\g'_k\oplus\dots\oplus\g'_n$. It follows that, for $x\in\mk{n}_n$, $Dx=nx+Z_{n+1}$, with $Z_{n+1}\in\mk{n}_{n+1}$, and similarly, $D'x=nx+Z'_{n+1}$ with $Z'_{n+1}\in\mk{n}_{n+1}$. Hence, $\left(D-D'\right)x=Z_{n+1}-Z'_{n+1}\in\mk{n}_{n+1}$. In other words, $D-D'\in\u\left(\mk{n}\right)$\ie $$D'\in D+\u\left(\mk{n}\right).$$ It follows, from Part \ref{orbitedeDsousU}), that there exists $X\in\u\left(\mk{n}\right)$ such that $D'=e^X\cdot D$. Then $$\mk{m}=\mk{n}^-(D')=\mk{n}^-(e^X\cdot D)=e^X.\mk{n}^-(D),$$ proving that the group $\U\left(\mk{n}\right)$ acts transitively on $\mk{n}^{\top}$.\\
Moreover, if $e^X.\mk{n}^-(D)=\mk{n}^-(D)$ with $X\in\u\left(\mk{n}\right)$, necessarly $e^X=\id$, and hence the action is simply transitive.

\end{enumerate}

\end{demo}

According to Theorem \ref{theoprincipal} and using the terminology of \cite{BL08}, $\left(\widetilde{\ca{F}},\widetilde{\ca{F}},\top\right)$ is a \emph{linear pair geometry}\ie for any filtration $\mk{n}\in\widetilde{\ca{F}}$, there is a transversal filtration $\mk{m}$, and the set $\mk{n}^\top$ carries a structure of $\K$-module by declaring the map $\u\left(\mk{n}\right) \rightarrow \mk{n}^\top,\ X \mapsto  e^X.\mk{m}$ to be a linear isomorphism.

\begin{coro}\label{coroprincipal}
If $\mk{n}=\mk{n}^+(\ad{E})\in\ca{F}$ is an inner filtration of $\g$, then we have $$\u\left(\mk{n}\right)=\ad{\mk{n}_1} \text{ and } \U\left(\mk{n}\right)=e^{\ad{\mk{n}_1}}.$$
Moreover, two transversal inner filtrations come from an inner grading of $\g$.
\end{coro}
\begin{demo}
For the first statement, we fix $X\in\u\left(\mk{n}\right)$ and prove by induction that, for all $2\pt n\pt 2k+1$, there are an element $Y^n\in\mk{n}_1$ and a derivation $R_n\in\u\left(\mk{n}\right)$ such that $\ad{Y^n}=X+R_n$ and $R_n(\mk{n}_i)\subset\mk{n}_{n+i}$, for all $-k+1\pt i\pt k$.\\
The second statement is proved in \cite{moi1}.
\end{demo}

\section{Generalized flag geometry}\label{varietededrapeaux}

\

In this section, we shall define a geometric object, called a \emph{generalized flag geometry}, associated with an inner $(2k+1)$-grading of $\g$, and realize this geometry by inner filtrations of $\g$.

\begin{defi}
We call a \emph{base point} of $\ca{F}\times\ca{F}$ a couple $(\mk{m},\mk{n})$ of transversal inner filtrations.\\ In this instance, according to Corollary \ref{coroprincipal}, $\mk{m}$ and $\mk{n}$ come from an inner grading of $\g$. In other words, $\mk{m}=\mk{n}^+(\ad{E})$ and $\mk{n}=\mk{n}^-(\ad{E})$. That is why, the base points of $\ca{F}\times\ca{F}$ will be denoted by $(\mk{n}^+, \mk{n}^-)$.
\end{defi}

\subsection{Projective elementary group}\label{groupprojelem}

Let $(\mk{n}^+,\mk{n}^-)$ be a base point of $\ca{F}\times\ca{F}$\ie we fix an inner grading of $\g$. Then $\mk{n}^+=\mk{n}^+(\ad{E})$ and $\mk{n}^-=\mk{n}^-(\ad{E})$. We consider the following subgroups of $\aut{\g}$: $$U^+:=\U\left(\mk{n}^+\right)=e^{\ad{\mk{n}^+_1}} \text{ and } U^-:=\U\left(\mk{n}^-\right)=e^{\ad{\mk{n}^-_1}}.$$
\begin{defi}
The group generated by $U^+$ and $U^-$, $G:=PE(\g_,E):=<U^+, U^->\subset \aut{\g}$, is called the \emph{projective elementary group} of $(\g,E)$.
\end{defi}
We also consider the following subgroups of $G$: the subgroup of $G$ preserving the grading, $$H:=\{g\in G, g\circ \ad{E}=\ad{E}\circ
g\},$$  and
$$P^+:=HU^+\ \text{ and }\
P^-:=HU^-.$$
Since the elements of $H$ commute with the characteristic derivation, they normalize $U^\pm$. Thus, $P^+$ and $P^-$ are indeed subgroups of $G$.\\

\subsection{Generalized flag geometry}

\

With notation as above, we define homogeneous spaces: $$M:=G/H,\quad
X^+:=G/P^- \text{ and } X^-:= G/P^+.$$ The base point $(P^+,P^-)$
of $X^-\times X^+$ is denoted by $(o^-,o^+)$. 

\begin{theo}\label{theostructure} With the above notation, we have:

\begin{enumerate}[(1)]

\item\label{H=,P=} The orbits of $D=\ad{E}$, respectively of
$\mk{n}^\pm$, under the action of $G$ are isomorphic to $M$,
respectively to $X^\mp$\ie
$$H=\{g\in G,\ g.(\mk{n}^+,\mk{n}^-)=(\mk{n}^+,\mk{n}^-)\}\quad
\text{and}\quad P^\pm=\{g\in G,\ g.\mk{n}^\pm=\mk{n}^\pm\}.$$

Moreover, $P^+\cap P^-=H$, $P^\pm\cap U^\mp=\{1\}$, and
$$P^\pm=\{g\in G,\ gDg^{-1}-D\in\ad{\mk{n}^\pm_1}\}.$$%

\item Let $x\in X^+$ and $y\in X^-$. If we identify $X^\pm$ with the corresponding orbits in
$\ca{F}$, then $x$ and $y$ are transversal if and only if there exists $g\in G$ such that $x=g.o^+$ and $y=g.o^-$.

\item For any $\mk{m}\in X^-$, we have $\mk{m}^\top\subset X^+$.
In particular, $\mk{n}^+_1$ is identified with $(o^-)^\top$ via $v\mapsto
e^{\ad{v}}.o^+$ so that we may consider that the space $\mk{n}^+_1$ is embedded in
$X^+$.

\item\label{omega+bijection} Let us consider the set
$\Omega^+:=\{g\in G,\ g.o^+\in \mk{n}^+_1\}$. The map $$\mk{n}^+_1\times
H\times \mk{n}^-_1 \rightarrow \Omega^+ ,\quad (v,h,w) \mapsto
e^{\ad{v}}he^{\ad{w}}$$ is a bijection.

\end{enumerate}

\end{theo}

\begin{demo}

\begin{enumerate}[(1)]
\item An element $g\in G$ stabilizes $(\mk{n}^+,\mk{n}^-)$ if and only if it
stabilizes the grading of $\g$, which means that $g$ commutes with
the derivation $D$\ie $g\in H$. Hence $$H=\{g\in G,\
g.(\mk{n}^+,\mk{n}^-)=(\mk{n}^+,\mk{n}^-)\}.$$

Moreover, $U^+$ and $H$ stabilize $\mk{n}^+$. This implies that  $P^+$ also stabilizes
$\mk{n}^+$. Conversely, if $g\in G$ satisfies $g.\mk{n}^+=\mk{n}^+$, then $g.\mk{n}^-$ is transversal to
$g.\mk{n}^+=\mk{n}^+$. It follows that there exists $v\in \mk{n}^+_1$ such that
$g.\mk{n}^-=e^{\ad{v}}.\mk{n}^-$. Hence $h:=e^{-\ad{v}}g$ stabilizes
$\mk{n}^+$ and $\mk{n}^-$\ie $h\in H$. Therefore $g=e^{\ad{v}}h\in
P^+$. The same also holds for $P^-$.

Thus $P^+\cap P^-$ is the stabilizer of $(\mk{n}^+,\mk{n}^-)$. It means that $P^+\cap P^-=H$.

Next, if $g\in P^+\cap U^-$, then $g=e^{\ad{v}}$, with $v\in \mk{n}^-_1$.
By Part \ref{Uagitsimplementtransitivement}) of Theorem
\ref{theoprincipal}, the map $v\mapsto e^{\ad{v}}.\mk{n}^+$ is
injective. Hence, since $g.\mk{n}^+=\mk{n}^+$, $v=0$ and $g=1$.

Finally, $g.\mk{n}^+=\mk{n}^+$ if and only if $D$ and $gDg^{-1}$
have the same plus-filtration, that is to say, from Section \ref{Transversalityandgradings}, if and only if
$gDg^{-1}-D\in\ad{\mk{n}^+_1}=\u\left(\mk{n}^+\right)$. Once again, the same applies for $P^-$.

\

\item Clearly, if $g\in G$, then $g.o^+\cong g.\mk{n}^-$ and $g.o^-\cong g.\mk{n}^+$ are transversal. Conversely, let $x\in X^+$ and $y\in X^-$ such that $x$ and $y$ are transversal. As $y\in X^-$, there exists $g\in G$ such that $y=g.o^-\cong g.\mk{n}^+$. Hence $g^{-1}y\cong \mk{n}^+$, and $g^{-1}x$ is transversal to $\mk{n}^+$. By Part \ref{Uagitsimplementtransitivement}) of
Theorem \ref{theoprincipal}, there exists $v\in \mk{n}^+_1$ such that $g^{-1}x=e^{\ad{v}}.\mk{n}^+$. Thus $x=ge^{\ad{v}}.o^+$, and $y=g.o^-=ge^{\ad{v}}.o^-$, since $e^{\ad{\mk{n}^+_1}}$ preserves $\mk{n}^+\cong o^-$.

\

\item Let $\mk{m}=g.o^-\in X^-$. By Part \ref{Uagitsimplementtransitivement}) of
Theorem \ref{theoprincipal}, $\mk{m}^\top\subset X^+$. In
particular, $o^-\in X^-$ and $(o^-)^\top=e^{\ad{\mk{n}^+_1}}.o^+$.

\

\item Let $(v,h,w)\in \mk{n}^+_1\times H\times \mk{n}^-_1$.We define
$g:=e^{\ad{v}}he^{\ad{w}}\in U^+P^-$. Hence $g.o^+\in \mk{n}^+_1$. It means that the
map is well-defined. Now, we assume that $g\in\Omega^+$ and we consider
$v:=g.o^+\in \mk{n}^+_1$. Then $e^{-\ad{v}}g.o^+=o^+$ and by Part
\ref{H=,P=}), $p:=e^{-\ad{v}}g\in P^-=HU^-$. This implies that there is $w\in \mk{n}^-_1$ such that $e^{-\ad{v}}g=he^{\ad{w}}$. Thus
$g=e^{\ad{v}}he^{\ad{w}}$, and the map is surjective.\\ Now, if $g=e^{\ad{v_1}}h_1e^{\ad{w_1}}=e^{\ad{v_2}}h_2e^{\ad{w_2}}$,
then $\underbrace{e^{-\ad{v_2}}e^{\ad{v_1}}h_1}_{\in
P^+}=\underbrace{h_2e^{\ad{w_2}}e^{-\ad{w_1}}}_{\in P^-}$. But, by
Part \ref{H=,P=}), $P^+\cap P^-=H$ so $v_1=v_2$. Similarly,
$w_1=w_2$ and it implies that $h_1=h_2$. Therefore, the map is injective too, hence is a bijection.

\end{enumerate}

\end{demo}

Actually, we have a $G$-equivariant imbedding $X^+\times X^-\hookrightarrow \ca{F}\times\ca{F}$.
\begin{defi}
The data $\left(X^+, X^-, \top\right)$ is called the \emph{generalized flag geometry (of $k$-graded type)} of $(\g,D)$.
\end{defi}

Moreover, according to Theorem \ref{theostructure}, $(X^+,X^-,\top)$ is also a linear pair geometry. Indeed, for any $x\in X^\pm$\exi $y\in X^\mp$ such that $x\top y$, and $x^\top$ carries a structure of $\K$-module. In particular, $$X^\pm=\bigcup_{y\in X^\mp}y^\top.$$ In other words, $$X^+=\bigcup_{g\in G} g\left(\mk{n}^+_1\right) \text{ and } X^-=\bigcup_{g\in G} g\left(\mk{n}^-_1\right).$$

\begin{defi}
We shall say that $$\ca{A}:=\left\{\left(g\left(\mk{n}^+_1\right), \varphi_g\right), g\in G\right\}, \text{ where } \varphi_g: g\left(\mk{n}^+_1\right)\subset X^+\rightarrow \mk{n}^+_1,\ g\cdot x\mapsto x,$$ is the \emph{natural atlas} of $X^+$ and the maps $\varphi_g$ will be called the \emph{charts} of the atlas $\ca{A}$.
\end{defi}
Now, let us define some operators which will be useful in the sequel.

\begin{defi}\label{defdeno}
For any $x\in \mk{n}^+_1$ and $g\in\aut{\g}$, we define the \emph{denominators}, resp. \emph{co-denominators}, of $g$ by:
$$d_g(x)_i=\text{pr}_{\mk{n}^+_i}\circ (e^{-\ad{x}}g^{-1}) \circ\iota_{\mk{n}^+_i}\in\text{End}(\mk{n}^+_i), \text{ resp. }
c_g(x)_i=\text{pr}_{\mk{n}^-_i}\circ (ge^{\ad{x}})
\circ\iota_{\mk{n}^-_i}\in\text{End}(\mk{n}^-_i), \text{ for } 1\pt i\pt k,$$ where
pr$_{\mk{n}^\pm_i}$ is the projection of $\g$ onto
$\mk{n}^\pm_i$ and $\iota_{\mk{n}^\pm_i}$ is the imbedding of
$\mk{n}^\pm_i$ into $\g$. Observe that, for any $1 \pt i\pt k$, the maps
$$d_g:\mk{n}^+_1\rightarrow\text{End}(\mk{n}^+_i),\ x \mapsto d_g(x)_i \text{ and }
c_g:\mk{n}^+_1\rightarrow\text{End}(\mk{n}^-_i),\ x\mapsto c_g(x)_i$$ are polynomial maps.\\
We also define, for $x\in \mk{n}^+_1$, $w\in \mk{n}^-_1$ and $1\pt i \pt k$,
$$B^+(x,w)_i:=d_{e^{\ad{w}}}(x)_i \text{ and }
B^-(w,x)_i:=c_{e^{\ad{w}}}(x)_i.$$ These operators are called
\emph{generalized Bergman operators}.
\end{defi}
Indeed, they can be considered as a generalization of Bergman operators of a Jordan pair \cite{Loos75}. More precisely, if $\g=TKK(V^+,V^-)$ is the Kantor-Koecher-Tits algebra of a Jordan pair $(V^+,V^-)$ (see \cite{Mey} for the definition of $TKK(V^+, V^-)$),then $k=1$ and the operators $B^+(x,w)=d_{e^{\ad{w}}}(x)$ and $B^-(w,x)=c_{e^{\ad{w}}}(x)$ are actually the Bergman operators of $(V^+,V^-)$ (see \cite{BN1}).
\begin{coro}\label{corog.xinV+}
Let $x\in \mk{n}^+_1$ and $g\in\aut{\g}$. Then the following statements are
equivalent:
\begin{enumerate}[i)]
\item $g\cdot x\in \mk{n}^+_1$\ie $\mk{n}^+$ and $ge^{\ad{x}}.\mk{n}^-$ are
transversal.
\item\label{concl} For all $i\in\{1,\dots,k\},\ d_g(x)_i \text{ and } c_g(x)_i$ are invertible.
\end{enumerate}
\end{coro}

\begin{demo}
By definition, $\mk{n}^+$ and $ge^{\ad{x}}.\mk{n}^-$ are transversal \ssi $\g=\mk{n}^+_n\oplus ge^{\ad{x}}\left(\mk{n}^-_{-n+1}\right)$, for all\\ $-k+1\pt n\pt k$. In other words, they are transversal if and only if, for all $1\pt i\pt k, ge^{\ad{x}}\left(\g_{-k}\oplus\dots\oplus\g_{-i}\right)$ is a complement of $\g_k\oplus\dots\oplus\g_{-i+1}$ and $\g_{-k}\oplus\dots\g_{i-1}$ is a complement of $e^{-\ad{x}}g^{-1}\left(\g_k\oplus\dots\oplus\g_i\right)$, and this last statement is equivalent to \ref{concl}).
\end{demo}

\

In particular, with the notation of Theorem \ref{theostructure},
we have $$\Omega^+=\{g\in G,\ \forall i\in\{1,\dots,k\},\
d_g(0)_i\text{ and } c_g(0)_i\text{ are
invertible}\}.$$

\section{Realization of $\g$ by polynomial vector fields}\label{g=polymap}

\

In this part, we define some algebraic objects that we call ``bundles'' on $\ca{F}$, even if we do not consider a topology. For example, we shall define the \emph{tangent} and \emph{structural bundle} of $\ca{F}$, as well as some sections of these bundles. Then, fixing a grading of $\g$, we shall consider the corresponding generalized flag geometry and construct a realization of $\g$ by polynomial fields on the positive part of the grading: $\mk{n}^+_1=\bigoplus\limits_{n=1}^k\g_n$. Finally, we shall give a global version of this realization\ie we shall study the action of the projective elementary group on $X^+$ in the ``chart'' $\mk{n}^+_1$.

\subsection{Bundles and vector fields}\label{bundles}

\begin{defi}\label{defTiX}
For any $(2k+1)$-filtration
$\mk{n}=(\mk{n}_k\subset\mk{n}_{k-1}\subset\dots\subset\mk{n}_{-k+1})\in\ca{F}$,
and any integer $1\pt i\pt k$, we define the following $\bb{K}$-modules:
$$T^{(i)}_{\mk{n}}\ca{F}:=\g/\mk{n}_{-i+1} \text{ and }
T'^{(i)}_{\mk{n}}\ca{F}:=\mk{n}_i.$$ To simplify the notations, instead of $T^{(1)}_{\mk{n}}\ca{F}$, we write $T_{\mk{n}}\ca{F}$ and call this module the \emph{tangent space of $\ca{F}$ at $\mk{n}$}, using the terminology of \cite{BN1}. Similarly, $T'^{(1)}_{\mk{n}}\ca{F}$ shall be denoted by $T'_{\mk{n}}\ca{F}$ and called the \emph{structural space of $\ca{F}$ at $\mk{n}$}.\\
Observe that, if $\mk{n}=\mk{n}^-(\ad{E})$ is the minus filtration defined by the Euler operator $E$, then\\ $\mk{n}_{-i+1}=\g_{-k}\oplus\dots\oplus\g_{i-1}$, $\mk{n}_{i}=\g_{-k}\oplus\dots\oplus\g_{-i}$,
and $$T^{(i)}_{\mk{n}}\ca{F}\cong\g_k\oplus\dots\oplus\g_i \text{ and }
T'^{(i)}_{\mk{n}}\ca{F}=\g_{-k}\oplus\dots\oplus\g_{-i}.$$
We  also define the disjoint unions:
$$T^{(i)}\ca{F}:=\bigcup_{\mk{n}\in\ca{F}}T^{(i)}_{\mk{n}}\ca{F} \text{ and } T'^{(i)}\ca{F}:=\bigcup_{\mk{n}\in\ca{F}}T'^{(i)}_{\mk{n}}\ca{F},$$ so that the following projections may be defined:
$$\pi^{(i)}: T^{(i)}\ca{F}  \rightarrow  \ca{F},\ Y\in T^{(i)}_{\mk{n}}\ca{F} \mapsto  \mk{n} \text{ and } \pi'^{(i)}:  T'^{(i)}\ca{F}  \rightarrow  \ca{F},\ x\in T'^{(i)}_{\mk{n}}\ca{F} \mapsto  \mk{n}.$$
We write $T\ca{F}:=T^{(1)}\ca{F}$ and
$T'\ca{F}:=T'^{(1)}\ca{F}$ and  call them the \emph{tangent bundle
of $\ca{F}$} and the \emph{structural bundle of $\ca{F}$}.
\end{defi}

\begin{lem}\label{cassemisimpledimfinie}
If $\g$ is a semi-simple, finite dimensional Lie algebra over a field $\bb{K}$, then, \prtt $\mk{n}\in\ca{F}$ and $1\pt i\pt k$, $\left(T^{(i)}_{\mk{n}}\ca{F}\right)^{\ast}\cong T'^{(i)}_{\mk{n}}\ca{F}$.
\end{lem}
\begin{demo}
Let $\mk{n}\in\ca{F}$ be an inner filtration. Since $\mk{n}$ is inner, it comes from a grading $\g=\g_k\oplus\dots\oplus\g_{-k}$. If $x\in\g_i$ and $y\in\g_j$, then $\ad{x}\ad{y}(\g_l)\subset\g_{i+j+l}$. But, $\g_l\cap\g_{i+j+l}=\{0\}$ if $i+j\neq 0$. It implies that, if we choose a basis of each $\g_n$, we obtain a basis of $\g$ in which the matrix of $\ad{x}\ad{y}$ has only zeros on the diagonal, so that the Killing form $B(x,y)=\tr{(\ad{x}\ad{y})}$ is zero. Hence, we can consider  $\widetilde{B}:\mk{n}_i\times\g/\mk{n}_{-i+1}\rightarrow \bb{K}$, which is non-degenerate because $\g$ is semi-simple. Therefore, $\left(T^{(i)}_{\mk{n}}\ca{F}\right)^\ast\cong T'^{(i)}_{\mk{n}}\ca{F}$.
\end{demo}
Thus, if $\g$ is a semi-simple, finite dimensional Lie algebra over a field, then the bundles $T'^{(i)}\ca{F}$ are dual bundles of $T^{(i)}\ca{F}$\ie $T'^{(i)}\ca{F}\cong\left(T^{(i)}\ca{F}\right)^\ast$. In particular, the structural bundle $T'\ca{F}$ is the \emph{cotangent bundle} $T^\ast\ca{F}$.

\

Since the group $\aut{\g}$ acts on $\ca{G}$ and $\ca{F}$, the following linear maps may be defined for any $g\in\aut{\g}$: $$\fct{$T^{(i)}_{\mk{n}}g$}{$T^{(i)}_{\mk{n}}\ca{F}$}{$T^{(i)}_{\mk{g.\mk{n}}}\ca{F}$}{$Y$
mod $\mk{n}_{-i+1}$}{$gY$ mod $g\mk{n}_{-i+1}$} \text{ and }
\fct{$T'^{(i)}_{\mk{n}}g$}{$T'^{(i)}_{\mk{n}}\ca{F}$}{$T'^{(i)}_{\mk{g.\mk{n}}}\ca{F}$}{$Y$}{$gY$.}$$
Then, defining the maps $T^{(i)}g$ and $T'^{(i)}g$ as usual, it holds:
$$T^{(i)}(g\circ h)=T^{(i)}g\circ
T^{(i)}h \text{ and } T'^{(i)}(g\circ h)=T'^{(i)}g\circ T'^{(i)}h.$$
Now, let us fix an inner grading of $\g$ and consider the associated generalized flag geometry. According to Theorem \ref{theostructure}, $X^\pm\subset\ca{F}$, so that we may define the spaces $T^{(i)}_{\mk{n}}X^\pm$, $T'^{(i)}_{\mk{n}}X^\pm$, as well as the bundles $T^{(i)}X^\pm$ and $T'^{(i)}X^\pm$. These bundles are, in an algebraic way, bundles associated to representations of $P^-$ given by the denominators at 0, $d_p(0)_i$, and  the co-denominators at 0, $c_p(0)_i$, introduced in Definition \ref{defdeno}.
More precisely, according to Corollary \ref{corog.xinV+}, for $g=p\in P^-$ and $x=0$, $d_p(0)_i$ and $c_p(0)_i$ are bijections and the following lemma holds. 
\begin{lem}
The maps $\rho^+_i: P^- \rightarrow \Gl(\mk{n}^+_i),\ p \mapsto d_p(0)_i^{-1}$ and $\rho^-_i: P^- \rightarrow \Gl(\mk{n}^-_i),\ p \mapsto c_p(0)_i$ are group homomorphisms.
\end{lem}
\begin{demo}
Let $p _1, p_2\in P^-$. We have $d_{p_1p_2}(0)_i=\text{pr}_{\mk{n}^+_i}\circ p_2^{-1}p_1^{-1}\circ\iota_{\mk{n}^+_i}$. Now, since $p_2\in P^-$, there is the identity: $$\text{pr}_{\mk{n}^+_i}\circ p_2^{-1}=\text{pr}_{\mk{n}^+_i}\circ p_2^{-1}\circ\text{pr}_{\mk{n}^+_i}.$$ We deduce that $d_{p_1p_2}(0)_i=\text{pr}_{\mk{n}^+_i}\circ p_2^{-1}\circ\text{pr}_{\mk{n}^+_i}\circ p_1^{-1}\iota_{\mk{n}^+_i}=d_{p_2}(0)_id_{p_1}(0)_i$. Therefore, $$d_{p_1p_2}(0)_i^{-1}=d_{p_1}(0)_i^{-1}d_{p_2}(0)_i^{-1}.$$\\
Moreover, it is obvious that $c_{p_1p_2}(0)_i=c_{p_1}(0)_ic_{p_2}(0)_i$.
\end{demo}

\

Now, let us focus our attention on $\rho_i^+$. We consider the quotient $$G\times_{P^-}\mk{n}^+_i:=G\times\mk{n}^+_i/\sim,$$ where $(g,x_i)\sim (gp,\rho^+_i(p)^{-1}x_i)$ for $g\in G, x_i\in\mk{n}^+_i$ and $p\in P^-$.
\begin{theo}\label{fibretangent}
For $1\pt i\pt k$, $T^{(i)}X^+\cong G\times_{P^-}\mk{n}^+_i$\ie there exists a $G$-equivariant bijection between $T^{(i)}X^+$ and $G\times_{P^-}\mk{n}^+_i$. We equally have $T'^{(i)}X^+\cong G\times_{P^-}\mk{n}^-_i$ and the same holds for $X^-$, if we exchange $P^-$ and $P^+$.
\end{theo}
\begin{demo}
In order to prove the identity $T^{(i)}X^+\cong G\times_{P^-}\mk{n}^+_i$, we consider the following map: $$ \varphi^{(i)}:  G\times_{P^-}\mk{n}^+_i \rightarrow T^{(i)}X^+,\ [g,x_i]  \mapsto  \left(T^{(i)}_{o^+}g\right)x_i=gx_i \text{ mod } g(\mk{n}^-_{-i+1}) $$
First, let us prove that $\varphi^{(i)}$ is well-defined. If $(gp,\rho^+_i(p)^{-1}x_i)$ is another representative of $[g,x_i]$, then $\left(T^{(i)}_{o^+}(gp)\right)(\rho^+_i(p)^{-1}x_i)=\left(T^{(i)}_{o^+}g\circ T^{(i)}_{o^+}p\right)d_p(0)_ix_i$. But the map $T^{(i)}_{o^+}p$ is defined by $$ T^{(i)}_{o^+}p:  T^{(i)}_{o^+}X^+\cong\mk{n}^+_i \rightarrow  T^{(i)}_{o^+}X^+\cong\mk{n}^+_i ,\ x_i  \mapsto  (\text{pr}_{\mk{n}^+_i}\circ p)x_i.$$ In other words, $T^{(i)}_{o^+}p=d_p(0)_i^{-1}$. It implies that $\left(T^{(i)}_{o^+}(gp)\right)\rho^+_i(p)^{-1}x_i=\left(T^{(i)}_{o^+}g\right)x_i$, so that $\varphi^{(i)}$ is well-defined.\\
Let $\mk{m}=g.o^+\in X^+$ and $x\in T^{(i)}_{\mk{m}}X^+=\g/\mk{m}_{-i+1}$. Since $g^{-1}x\in T^{(i)}_{o^+}X^+\cong\mk{n}^+_i$, it follows, from $\varphi^{(i)}\left([g,g^{-1}x]\right)=x\in T^{(i)}_{\mk{m}}X^+$, that $\varphi^{(i)}$ is surjective.\\
Now, let $(g_1,g_2,x_1,x_2)$ such that $\left(T^{(i)}_{o^+}g_1\right)x_1=\left(T^{(i)}_{o^+}g_2\right)x_2$.\\
First, we have $g_1.o^+=g_2.o^+$\ie $g_1^{-1}g_2.o^+=o^+$ so $g_1^{-1}g_2=p\in P^-$\ie $g_2=g_1p$. Hence\\ $\left(T^{(i)}_{o^+}g_2\right)x_2=\left(T^{(i)}_{o^+}g_1\circ T^{(i)}_{o^+}p\right)x_2=\left(T^{(i)}_{o^+}g_1\right)d_p(0)_i^{-1}x_2$. Therefore $x_1-d_p(0)_i^{-1}x_2\in\text{Ker}\left(T^{(i)}_{o^+}g_1\right)$. But, if $x\in\mk{n}^+_i$ is in the kernel of $T^{(i)}_{o^+}g_1$, then $g_1x\in g(\mk{n}^-_{-i+1})$\ie $x\in\mk{n}^+_i\cap\mk{n}^-_{-i+1}=\{0\}$. Thus, $x_1=d_p(0)_i^{-1}x_2$ and $(g_1,x_1)\sim(g_2,x_2)$, so that $\varphi^{(i)}$ is injective.\\
Finally, $\varphi^{(i)}\left(g'\cdot[g,x_i]\right)=\varphi^{(i)}\left([g'g,x_i]\right)=\left(T^{(i)}_{g.o^+}g'\circ T^{(i)}_{o^+}g\right)x_i=\left(T^{(i)}_{g.o^+}g'\right)\cdot\varphi^{(i)}\left([g,x_i]\right)$. Thus $\varphi^{(i)}$ is $G$-equivariant, and $T^{(i)}X^+$ and $G\times_{P^-}\mk{n}^+_i$ are isomorphic $G$-bundles.\\
The proofs are similar for the other statements.
\end{demo}

\

Now, let us focus on sections on these bundles. In fact, sections of $T^{(i)}\ca{F}$ can be easily constructed. Indeed, the following holds:
\begin{defi}\label{defsectionY}
For any $Y\in\g$ and $1\pt i\pt k$, we define the map
$$ \widetilde{Y}^{(i)}:  \ca{F}  \rightarrow  T^{(i)}\ca{F},\ \mk{n}  \mapsto  Y^{(i)}_\mk{n},$$ where $Y^{(i)}_\mk{n}:= Y\ \text{mod } \mk{n}_{-i+1} \in
T^{(i)}_{\mk{n}}\ca{F}.$\\
In particular, for $i=1$, $Y_\mk{n}:=Y^{(1)}_\mk{n}$ is called the \emph{value of $Y$ at $\mk{n}$}, and the map $ \widetilde{Y}: \ca{F} \rightarrow T\ca{F},\ \mk{n}  \mapsto Y_{\mk{n}}$
defines a \emph{vector field on $\ca{F}$}, called a \emph{projective vector field}.
\end{defi}
\begin{lem}
The maps $\widetilde{Y}^{(i)}$ are sections of $T^{(i)}\ca{F}$.
\end{lem}
\begin{demo}
By definition, $\pi^{(i)}\circ \widetilde{Y}^{(i)}(\mk{n})=\pi^{(i)}(Y\text{ mod } \mk{n}_{-i+1})=\mk{n}$. Hence $\pi^{(i)}\circ\widetilde{Y}^{(i)}=\text{Id}_{\ca{F}}$.
\end{demo}

\subsection{The canonical kernel}\label{noyaucanonique}

\

Let $(\mk{m},\mk{n})\in\ca{F}$. The spaces $T_{\mk{n}}\ca{F}$ and $T'_{\mk{m}}\ca{F}$ are filtered $\bb{K}$-modules:
$$\begin{array}{ccccccccccc} T_{\mk{n}}\ca{F}=T^{(1)}_{\mk{n}}\ca{F} & \twoheadrightarrow & T^{(2)}_{\mk{n}}\ca{F} & \twoheadrightarrow & \dots & \twoheadrightarrow & T^{(i)}_{\mk{n}}\ca{F} & \twoheadrightarrow & \dots & \twoheadrightarrow & T^{(k)}_{\mk{n}}\ca{F} \quad \text{ and }\\ \\ T'^{(k)}_{\mk{m}}\ca{F} & \hookrightarrow & T'^{(k-1)}_{\mk{m}}\ca{F} & \hookrightarrow & \dots & \hookrightarrow & T'^{(k-i)}_{\mk{m}}\ca{F} & \hookrightarrow & \dots & \hookrightarrow & T'^{(1)}_{\mk{m}}\ca{F}=T'_{\mk{m}}\ca{F}.\\ \end{array}$$
Observe that, for $1\pt i\pt k$, $T'^{(i)}_{\mk{m}}\ca{F}$ is a subspace of $\g$ which defines sections of $T^{(i)}\ca{F}$, so that we may define the following $\bb{K}$-linear maps:
$$ K^{(i)}_{\mk{n},\mk{m}}:  T'^{(i)}_{\mk{m}}\ca{F}=\mk{m}_i  \rightarrow  T^{(i)}_{\mk{n}}\ca{F}=\mk{g}/\mk{n}_{-i+1},\
 Y  \mapsto  Y^{(i)}_{\mk{n}}= Y \text{mod } \mk{n}_{-i+1}.$$

By definition, there is the following diagram:

$$\xymatrix{T'^{(1)}_{\mk{m}}\ca{F}=\mk{m}_1\ar^{\hspace{-.5cm}K^{(1)}_{\mk{n},\mk{m}}}[r] & \g/\mk{n}_0=T^{(1)}_{\mk{n}}\ca{F} \ar@{->>}[d]\\
T'^{(2)}_{\mk{m}}\ca{F}=\mk{m}_2 \ar@{^(->}[u]
\ar^{\hspace{-.5cm}K^{(2)}_{\mk{n},\mk{m}}}[r] & \g/\mk{n}_{-1}=T^{(2)}_{\mk{n}}\ca{F} \ar@{->>}[d]\\
 \vdots \ar@{^(->}[u] & \vdots \ar@{->>}[d] \\
T'^{(i)}_{\mk{m}}\ca{F}=\mk{m}_i
\ar_{\hspace{-.6cm}K^{(i)}_{\mk{n},\mk{m}}}[r]
\ar@{^(->}[u] & \g/\mk{n}_{-i+1}=T^{(i)}_{\mk{n}}\ca{F} \ar@{->>}[d]\\
\vdots \ar@{^(->}[u] & \vdots \ar@{->>}[d] \\
T'^{(k)}_{\mk{m}}\ca{F}=\mk{m}_k \ar_{\hspace{-.6cm}
K^{(k)}_{\mk{n},\mk{m}}}[r] \ar@{^(->}[u] & \g/\mk{n}_{-k+1}=T^{(k)}_{\mk{n}}\ca{F}\\}$$
\begin{defi}\label{noyaucano}
The collection of maps
$\left(K^{(i)}_{\mk{m},\mk{n}},K^{(i)}_{\mk{n},\mk{m}}\right)$, for
$(\mk{m}, \mk{n})\in\ca{F}\times\ca{F}$, and $1\pt i\pt
k$, is called the \emph{canonical kernel}.
\end{defi}
We shall now prove that $K^{(i)}:(\mk{m},\mk{n})\mapsto K^{(i)}_{\mk{n},\mk{m}}$ defines a section of a bundle over $\ca{F}\times\ca{F}$. To this end, we define \begin{center}$\hom{T'\ca{F},T\ca{F}}:=\displaystyle{\bigcup_{\mk{m},\mk{n}\in\ca{F}}}\Hom\nolimits_{\bb{K}}\left(T'^{(i)}_{\mk{m}}\ca{F},T^{(i)}_{\mk{n}}\ca{F}\right),$\end{center} where the union is disjoint, and the projection
\begin{center}$\Pi^{(i)}:  \hom{T'\ca{F},T\ca{F}}  \rightarrow  \ca{F}\times\ca{F},\ \varphi_{\mk{n},\mk{m}}\in\Hom_{\bb{K}}\left(T'^{(i)}_{\mk{m}}\ca{F},T^{(i)}_{\mk{n}}\ca{F}\right)  \mapsto  (\mk{m},\mk{n}).$\end{center}
By definition, $K^{(i)}_{\mk{n},\mk{m}}\in\Hom\left(T'^{(i)}_{\mk{m}}\ca{F},T^{(i)}_{\mk{n}}\ca{F}\right)$ and we define \begin{center}$K^{(i)}:  \ca{F}\times\ca{F} \rightarrow \hom{T'\ca{F},T\ca{F}},\ (\mk{m},\mk{n}) \mapsto K^{(i)}_{\mk{n},\mk{m}}.$\end{center} Thus $\Pi^{(i)}\circ K^{(i)}=\id_{\ca{F}\times\ca{F}}$.
Hence, $K^{(i)}$ is a section of the bundle \raisebox{.55cm}{$\xymatrix{\hom{T'\ca{F},T\ca{F}} \ar@{->>}[d] \\ \ca{F}\times\ca{F}}$}.\\[.2cm]
\begin{rem}
According to Lemma \ref{cassemisimpledimfinie}, if $\g$ is a semi-simple, finite dimensional Lie algebra over a field, the space $T'^{(i)}_{\mk{m}}\ca{F}$ is isomorphic to $\left(T^{(i)}_{\mk{m}}\ca{F}\right)^{\ast}$. Hence \begin{center} $\hom{T'\ca{F},T\ca{F}}\cong\left(T'^{(i)}\ca{F}\right)^\ast\boxtimes T^{(i)}\ca{F}=T^{(i)}\ca{F}\boxtimes T^{(i)}\ca{F}.$\end{center}
In other words, $\hom{T'\ca{F},T\ca{F}}$ is a bundle over $\ca{F}\times\ca{F}$ and if $(\mk{m},\mk{n})\in\ca{F}$ are two filtrations, the fiber over $(\mk{m},\mk{n})$ is isomorphic to  $T^{(i)}_{\mk{m}}\ca{F}\otimes T^{(i)}_{\mk{n}}\ca{F}$.\\
\end{rem}
In general case, the following holds:
\begin{theo}\label{noyaucaniosection}
 \begin{enumerate}[(1)]
  \item The map $K^{(i)}$ is a $G$-equivariant section of $\hom{T'\ca{F},T\ca{F}}$.
  \item Let $(\mk{m},\mk{n})\in\ca{F}\times\ca{F}$ be two inner filtrations. Then, they are transversal if and only if, \prtt $i=1,\dots, k$, the maps $K^{(i)}_{\mk{m},\mk{n}}$ and $K^{(i)}_{\mk{n},\mk{m}}$ are isomorphisms.
 \end{enumerate}
\end{theo}
\begin{demo}
\begin{enumerate}[(1)]
\item The projective elementary group $G$ acts on $\hom{T'\ca{F},T\ca{F}}$ as follows: if $g\in G$ and $\varphi\in\Hom\left(T'^{(i)}_{\mk{m}}\ca{F},T^{(i)}_{\mk{n}}\ca{F}\right)$, then $$g\cdot\varphi=T^{(i)}_{\mk{n}}g\circ\varphi\circ\left(T'^{(i)}_{\mk{m}}g\right)^{-1}\in\Hom\left(T'^{(i)}_{g.\mk{m}}\ca{F},T^{(i)}_{g.\mk{n}}\ca{F}\right).$$ Since we have already noted that $K^{(i)}$ is a section, we only have to prove that it is equivariant. Let $(\mk{m},\mk{n})\in\ca{F}\times\ca{F}$ and $g\in G$. The map $K^{(i)}_{g.\mk{n},g.\mk{m}}$ is defined by:
$$K^{(i)}_{g.\mk{n},g.\mk{m}}:  T'^{(i)}_{g.\mk{m}}\ca{F}  \rightarrow  T^{(i)}_{\mk{n}}\ca{F},\ Y\in g(\mk{m}_i)  \mapsto  Y \text{ mod } g(\mk{n}_{-i+1}).$$
Let $Y\in g(\mk{m}_i)$. There exists $Z\in\mk{m}_i$ such that $Y=gZ=\left(T'^{(i)}_{\mk{m}}g\right)Z$, and\\
$\begin{array}{ccc} K^{(i)}_{g.\mk{n},g.\mk{m}}(gZ) & = & gZ \text{ mod } g(\mk{n}_{-i+1})\\ & = & \left(T^{(i)}_{\mk{n}}g\right)\left(Z \text{ mod } \mk{n}_{-i+1}\right)\\ & = & \left(T^{(i)}_{\mk{n}}g\right)\left(K^{(i)}_{\mk{n},\mk{m}}(Z)\right)\\ \end{array}$\\
In other words, $$K^{(i)}_{g.\mk{n},g.\mk{m}}=T^{(i)}_{\mk{n}}g\circ K^{(i)}_{\mk{n},\mk{m}}\circ \left(T'^{(i)}_{\mk{m}}g\right)^{-1}=g\cdot K^{(i)}_{\mk{n},\mk{m}},$$\textit{i.e.} $K^{(i)}$ is $G$-equivariant.
\item The map $K^{(i)}_{\mk{n},\mk{m}}$ is bijective \ssi $\g=\mk{m}_i\oplus\mk{n}_{-i+1}$. Thus $\mk{m}$ and $\mk{n}$ are transversal if and only if, for all $1\pt i\pt k$, the maps $K^{(i)}_{\mk{n},\mk{m}}$ and $K^{(i)}_{\mk{m},\mk{n}}$ are isomorphisms.
\end{enumerate}
\end{demo}

From Theorem \ref{theostructure}, there is an imbedding $X^\pm \subset \ca{F}$. Therefore, we can consider $\hom{T' X^+, TX^-}$ or $\hom{T'X^-,TX^+}$. These bundles are isomorphic to ``associated bundles'' (like $T^{(i)}X^\pm$ in Theorem \ref{fibretangent}). Indeed, the following holds:
\begin{theo}
For all $i=1,\dots, k$, we have $${\rm Hom}^{(i)}\left(T'X^+,TX^-\right)\cong (G\times G)\times_{\left(P^-\times P^+\right)}\End{\mk{n}^-_i} \text{ and } {\rm Hom}^{(i)}\left(T'X^-,TX^+\right)\cong(G\times G)\times_{\left(P^+\times P^-\right)}\End{\mk{n}^+_i},$$
where $(G\times G)\times_{\left(P^-\times P^+\right)}\End{\mk{n}^-_i}:=(G\times G)\times \End{\mk {n}^-_i}/\sim$ with $(g, g', \varphi)\sim (gp^-, g'p^+, c_{p^+}(0)_i^{-1}\varphi c_{p^-}(0)_i)$\\ and $(G\times G)\times_{\left(P^+\times P^-\right)}\End{\mk{n}^+_i}:=(G\times G)\times \End{\mk {n}^+_i}/\sim$ with $(g, g', \varphi)\sim (gp^+, g'p^-, d_{p^+}(0)_i\varphi d_{p^-}(0)_i^{-1})$.
\end{theo}

\begin{demo}
The proof is similar to the one of Theorem \ref{fibretangent}.
\end{demo}

Now we can consider restrictions of canonical kernel $K^{(i)}$ over $X^\pm\times X^\mp$ and regard it as a section of $\hom{T'X^\pm,TX^\mp}$. Moreover, if $x\in\mk{n}^+_1\subset X^+$ and $y\in\mk{n}^-_1\subset X^-$, the maps $K^{(i)}_{x,y}$ are given by $$K^{(i)}_{x,y}:  T'^{(i)}_{y\cdot o^-}X^-  \rightarrow  T^{(i)}_{x\cdot o^+}X^+,\ e^{\ad{y}}Y  \mapsto  \text{pr}_{\mk{n}^+_1}\left(e^{-\ad{x}}e^{\ad{y}}Y\right).$$ In other words, $K^{(i)}_{x,y}=B^+(x,-y)_1$\ie $K^{(i)}_{x,y}$ is a generalized Bergman operator.

\subsection{Realization of $\g$ by polynomial fields}\label{champspoly}

\

Now, let us fix an inner $(2k+1)$-grading of $\g$ and consider the general flag geometry associated with $(\g,D)$. In the previous part, we proved that $T^{(i)}X^+$ is isomorphic to $G\times_{P^-}\mk{n}^+_i$. It allows us to identify $T^{(i)}X^+$ and $G\times_{P^-}\mk{n}^+_i$. In addition, the maps $\widetilde{Y}^{(i)}$ define sections of $T^{(i)}X^+$. But, we know that the sections of $G\times_{P^-}\mk{n}^+_i$ are identified with the induced representation of $G$, so that they correspond to the functions $f:G\rightarrow \mk{n}^+_i$ such that, \prtt $p\in P^-$, $f(gp)=d_p(0)_if(g)$. Indeed, let us consider such a function $f$ and define $$s_f: X^+  \rightarrow  T^{(i)}X^+,\ g.o^+  \mapsto  \left(T^{(i)}_{o^+}g\right)f(g).$$ Then $s_f$ is a section of \fibre{T^{(i)}X^+}{X^+}, and all sections arise in this way. Precisely, the function corresponding to $\widetilde{Y}^{(i)}$ is $$f_Y^{(i)}:  G  \rightarrow  \mk{n}^+_i,\  g \mapsto  \text{pr}_{\mk{n}^+_i}\left(g^{-1}Y\right).$$
In particular, if $g=e^{\ad{x}}$, with $x\in\mk{n}^+_i$, we identify the restriction of $\widetilde{Y}^{(i)}$ on $\mk{n}^+_i$ with $$\widetilde{Y}^{+(i)}=f_{Y|_{\mk{n}^+_i}}: \mk{n}^+_i \rightarrow \mk{n}^+_i,\ x  \mapsto  \text{pr}_{\mk{n}^+_i}\left(e^{-\ad{x}}Y\right).$$
We remark that $x\mapsto \widetilde{Y}^{+(i)}(x)$ is a polynomial map of $x$. Thus, we obtain a map $$\g  \rightarrow  Pol\left(\mk{n}^+_i\right),\ Y  \mapsto  \widetilde{Y}^{+(i)},$$ where $Pol\left(\mk{n}^+_i\right)$ is the space of polynomial self-mappings of $\mk{n}^+_i$. In particular for $i=1$, we obtain a realization of $\g$ by polynomial fields on $\mk{n}^+_1$.
For example, if $\g$ is 3-graded, there is only one such map:\\ $\widetilde{Y}^+=\widetilde{Y}^{+(1)}: x\mapsto \text{pr}_{\g_1}\left(e^{-\ad{x}}Y\right)=\text{pr}_{\g_1}\left(Y-\left[x,Y\right]+\dfrac{1}{2}\left[x,\left[x,Y\right]\right]\right)$ and we recover formulas of \cite{BN1}.\\[.2cm]
If $\g$ is 5-graded, there are two such maps:\\
$\widetilde{Y}^{+(2)}: x\in\g_2\mapsto \text{pr}_{\g_2}\left(e^{-\ad{x}}Y\right)=\text{pr}_{\g_2}\left(Y-\left[x,Y\right]+\dfrac{1}{2}\left[x,\left[x,Y\right]\right]\right)$. Hence, we obtain the formulas corresponding to the 3-grading $\g_2\oplus\g_0\oplus\g_2$:\\
$\begin{array}{l}
\text{if } Y\in\g_2, \widetilde{Y}^{+(2)}(x)=Y,\\
\text{if } Y\in\g_1, \widetilde{Y}^{+(2)}(x)=0,\\
\text{if } Y\in\g_0, \widetilde{Y}^{+(2)}(x)=-\left[x,Y\right],\\
\text{if } Y\in\g_{-1}, \widetilde{Y}^{+(2)}(x)=0,\\
\text{if } Y\in\g_{-2}, \widetilde{Y}^{+(2)}(x)=\dfrac{1}{2}\left[x,\left[x,Y\right]\right].\\
\end{array}$\\
Note that we also have $\widetilde{Y}^{+(1)}: x=x_1+x_2\in\g_1\oplus\g_2\mapsto \text{pr}_{\g_1\oplus\g_2}\left(e^{-\ad{x}}Y\right)$\ie\\ $\widetilde{Y}^{+(1)}(x)=\text{pr}_{\g_1\oplus\g_2}\left(Y-\left[x,Y\right]+\dfrac{1}{2}\left[x,\left[x,Y\right]\right]-\dfrac{1}{6}\left[x,\left[x,\left[x,Y\right]\right]\right]+\dfrac{1}{24}\left[x,\left[x,\left[x,\left[x,Y\right]\right]\right]\right]\right).$\\
More precisely,\\ $\begin{array}{l}
\text{if } Y\in\g_2, \widetilde{Y}^{+(1)}(x)=Y,\\
\text{if } Y\in\g_1, \widetilde{Y}^{+(1)}(x)=Y-\left[x,Y\right],\\
\text{if } Y\in\g_0, \widetilde{Y}^{+(1)}(x)=-\left[x_1,Y\right]-\left[x_2,Y\right]+\dfrac{1}{2}\left[x_1,\left[x_1,Y\right]\right],\\
\text{if } Y\in\g_{-1}, \widetilde{Y}^{+(1)}(x)=-\left[x_2,Y\right]+\dfrac{1}{2}\left[x_1,\left[x_1,Y\right]\right]+\left[x_1,\left[x_2,Y\right]\right]-\dfrac{1}{6}\left[x_1,\left[x_1,\left[x_1,Y\right]\right]\right],\\
\text{if } Y\in\g_{-2}, \widetilde{Y}^{+(2)}(x)=\raisebox{-.35cm}{$\begin{array}{l} \left[x_1,\left[x_2,Y\right]\right]-\dfrac{1}{6}\left[x_1,\left[x_1,\left[x_1,Y\right]\right]\right]+\\ \dfrac{1}{2}\left[x_2,\left[x_2,Y\right]\right]-\dfrac{1}{2}\left[x_1,\left[x_1,\left[x_2,Y\right]\right]\right]+\dfrac{1}{24}\left[x_1,\left[x_1,\left[x_1,\left[x_1,Y\right]\right]\right]\right].\\ \end{array}$}\\
\end{array}$\\[.2cm]

\subsection{Description of the $G$-action in the chart $\mk{n}^+_1$}\label{actiondeG}

\subsubsection{Example: action by homographies}\label{actionparhomographie}

\

First, let us consider the example of $\g=\mk{gl}_n(\C)$ with the 3-grading defined by: $$\g_1=\left\{\left(\begin{array}{cc} 0 & B \\ 0 & 0 \\
\end{array}\right), B\in\ca{M}_{p,n-p}(\C)\right\},\
\g_0=\left\{\left(\begin{array}{cc} A & 0 \\ 0 & D \\
\end{array}\right), A\in\ca{M}_p(\C), D\in\ca{M}_{n-p}(\C)\right\},$$
$$\text{and } \g_{-1}=\left\{\left(\begin{array}{cc} 0 & 0 \\ C & 0
\\ \end{array}\right), C\in\ca{M}_{n-p,p}(\C)\right\}.$$
In this case, $G=\text{Gl}_n(\C)$ and $$X^-\cong Gras_p(\C^n) \text{ and } X^+\cong Gras_{n-p}(\C^n),$$\textit{i.e.} $X^-$ is the set of $p$-dimensional subspaces of $\C^n$ and $X^+$ the one of $n-p$ dimensional subspaces. In addition, $\mk{n}^+_1\cong\ca{M}_{p,n-p}(\C)$ and $\mk{n}^+_1\hookrightarrow X^+$ via $X\mapsto \left\{\left(\begin{array}{c} v\\ Xv\\ \end{array}\right), v\in\C^{n-p}\right\}$. In this chart, the action of the group $G=\Gl_n(\C)$ is given by $$\left(\begin{array}{cc} A & B \\ C & D \\ \end{array}\right)\cdot\left(\begin{array}{c} v \\ Xv \\ \end{array}\right)=\left(\begin{array}{c} v \\ (C+DX)(A+BX)^{-1}v\\ \end{array}\right),$$ if $(A+BX)$ is invertible. This condition means that the point $x=\left(\begin{array}{c} v \\ Xv \\ \end{array}\right)\in Gras_p(\C^n)$ is not sent to infinity by $g=\left(\begin{array}{cc} A & B \\ C & D \\ \end{array}\right)$. In other words, the projective group acts, in the chart $\mk{n}^+_1$, by homographies.

\subsubsection{General case: action by ``birational'' maps}

\

In the sequel, we shall focus our attention on the action of the projective elementary group $G$ in the chart $\mk{n}^+_1$. More precisely, for $g\in\aut{\g}$ and $x\in\mk{n}^+_1$ such that $g\cdot x\in\mk{n}^+_1$, $g\cdot x$ has a ``birational'' expression. For example, if $\g$ is 3-graded, by \cite{BN1}, it holds: $$g\cdot x=d_g(x)^{-1}n_g(x), \text{ where } x \mapsto n_g(x) \text{ is a quadratic polynomial map.}$$ When the grading is longer, the following holds:
\begin{prop}\label{propcocycle}
Let $g\in \aut{\g}$, and $x\in\mk{n}^+_i$ such that $g\cdot x\in\mk{n}^+_i$.
Then
$$\widetilde{\left(g^{-1}Y\right)}^{+(i)}(x)=d_g(x)_ie^{\ad{g\cdot x}}\widetilde{Y}^{+(i)}(g\cdot x).$$
Moreover, if $x\in\mk{n}^+_i$, and $g_1, g_2\in \aut{\g}$ are such that $g_2 \cdot x\in\mk{n}^+_i$ and $g_1g_2\cdot x\in\mk{n}^+_i$, then
$$d_{g_1g_2}(x)_i=d_{g_2}(x)_ie^{\ad{g_2\cdot x}} d_{g_1}(g_2\cdot x)_i=d_{g_2}(x)_id_{\text{Id}}(g_2\cdot x)_i^{-1}d_{g_1}(g_2\cdot x)_i.$$
\end{prop}
\begin{demo}
Let $x\in\mk{n}^+_i$ and $g\in\aut{\g}$ such that $g\cdot
x\in\mk{n}^+_i$. By Part \ref{omega+bijection}) of Theorem \ref{theostructure}\exi $p(g,x)\in P^-$ such that
$ge^{\ad{x}}=e^{\ad{g\cdot x}}p(g,x)$. It follows that $$p(g,x)^{-1}=e^{-\ad{x}}g^{-1}e^{\ad{g\cdot x}},$$ which implies that
$$p(g,x)^{-1}_i:=\text{pr}_{\mk{n}^+_i}\circ p(g,x)^{-1}\circ\iota_{\mk{n}^+_i}=d_g(x)_ie^{\ad{g\cdot
x}}\circ\iota_{\mk{n}^+_i}=d_g(x)_id_{\text{Id}}(g\cdot x)_i^{-1}$$ because
$g\cdot x\in\mk{n}^+_1$. Hence,
$$\widetilde{\left(g^{-1}Y\right)}^{+(i)}(x)=\text{pr}_{\mk{n}^+_i}\left(e^{-\ad{x}}g^{-1}Y\right)=\text{pr}_{\mk{n}^+_i}\left(p(g,x)^{-1}e^{-\ad{g\cdot
x}}Y\right).$$
Moreover, since $p(g,x)^{-1}\in P^-$, we have $$\text{pr}_{\mk{n}_i^+}\circ
p(g,x)^{-1}=p(g,x)_i^{-1}\circ \text{pr}_{\mk{n}^+_i},$$\textit{i.e.} for
$y\in\g$, the components in $\mk{n}^+_i$ of $p(g,x)^{-1}y$ depend only on 
components of $y$ in $\mk{n}^+_i$. Thus
$$\widetilde{\left(g^{-1}Y\right)}^{+(i)}(x)=d_g(x)_ie^{\ad{g\cdot
x}}\text{pr}_{\mk{n}^+_i}\left(e^{-\ad{g\cdot
x}}Y\right)=d_g(x)_ie^{\ad{g\cdot
x}}\widetilde{Y}^{+(i)}(g\cdot x).$$
In particular,
for $Y=v\in\mk{n}^+_i$, we have
$$\widetilde{\left(g^{-1}v\right)}^{+(i)}(x)=d_g(x)_ie^{\ad{g\cdot
x}}\text{pr}_{\mk{n}^+_i}\left(e^{-\ad{g\cdot
x}}v\right)=d_g(x)_iv.$$
Now, if $g_2\cdot x\in\mk{n}^+_i$ and $g_1g_2\cdot
x\in\mk{n}^+_i$, then
$$\begin{array}{rcl}
d_{g_1g_2}(x)_iv=\widetilde{\left(g_2^{-1}g_1^{-1}v\right)}^{+(i)}(x) & = &
\displaystyle{d_{g_2}(x)_ie^{\ad{g_2\cdot x}}\widetilde{\left(g_1^{-1}v\right)}^{+(i)}(g_2\cdot x)}\\
& = & \displaystyle{d_{g_2}(x)_ie^{\ad{g_2\cdot x}}d_{g_1}(g_2\cdot x)_iv.}\\ \end{array}$$
Hence $\displaystyle{d_{g_1g_2}(x)_i=d_{g_2}(x)_ie^{\ad{g_2\cdot x}}d_{g_1}(g_2\cdot x)_i}$.
\end{demo}
Now, we can give a ``rational'' expression for the action of $G$ on $\mk{n}^+_1\subset X^+$.
\begin{defi}
Let $g\in\aut{\g}$ and $x\in\mk{n}^+_1$. We define
$$n_g(x):=\text{pr}_{\mk{n}^+_1}\left(e^{-\ad{x}}g^{-1}E\right),$$ where $E$ is the Euler operator of the grading of $\g$. Observe that
$n_g(x)=\widetilde{\left(g^{-1}E\right)}^{+(1)}(x)$ and the map $n_g:\mk{n}^+_1\rightarrow \mk{n}^+_1,\ x\mapsto
n_g(x)$, called the \emph{numerator of $g$}, is polynomial.
\end{defi}
\begin{prop}\label{actionrationnelle}
Let  $g\in\aut{\g}$ and $x\in\mk{n}^+_1$ such that $g\cdot x\in\mk{n}^+_1$. Then $$E-e^{\ad{g\cdot x}}E=d_g(x)_1^{-1}n_g(x).$$
More precisely, there is a map $\psi:\mk{n}^+_1\rightarrow \mk{n}^+_1$ such thaht $$[E,g\cdot x]=d_g(x)_1^{-1}n_g(x)+\psi(g\cdot x).$$
\end{prop}
\begin{demo}
By definition, $\widetilde{\left(g^{-1}E\right)}^{+(1)}(x)=n_g(x)$, and by Proposition \ref{propcocycle}, $$\widetilde{\left(g^{-1}E\right)}^{+(1)}(x)=d_g(x)_1e^{\ad{g\cdot x}}\text{pr}_{\mk{n}^+_1}\left(e^{-\ad{g\cdot x}}E\right).$$
In other words, $$d_g(x)_1e^{\ad{g\cdot x}}\text{pr}_{\mk{n}^+_1}\left(e^{-\ad{g\cdot x}}E\right)=n_g(x),$$\textit{i.e.} $e^{\ad{g\cdot x}}\text{pr}_{\mk{n}^+_1}\left(e^{-\ad{g\cdot x}}E\right)=d_g(x)_1^{-1}n_g(x)$.
Since $g\cdot x\in\mk{n}^+_1$, pr$_{\mk{n}^+_1}\left(e^{-\ad{g\cdot x}}E\right)=e^{-\ad{g\cdot x}}E-E$ which implies that $$e^{\ad{g\cdot x}}\text{pr}_{\mk{n}^+_1}\left(e^{-\ad{g\cdot x}}E\right)=E-e^{\ad{g\cdot x}}E.$$ Hence $$E-e^{\ad{g\cdot x}}E=d_g(x)^{-1}_1n_g(x).$$
Now, let us fix $v\in\mk{n}^+_1$ and write $v=\sum_{n=1}^kv_n$, with $v_n\in\g_n$. Then $$E-e^{\ad{v}}E=[E,v]+\frac{1}{2}\left[v,\left[E,v\right]\right]+\dots=\sum_{n=1}^knv_n+\psi_n(v_1,\dots,v_{n-1}),$$ where the map $\psi_n:\g_1\times\g_2\times\dots\times\g_{n-1}\rightarrow \g_n$ is $(n-1)$-linear and is the part of the component of $E-e^{\ad{v}}E$ in $\g_n$ which only depends on $v_1,\dots,v_{n-1}$. For example, we have $\psi_1,\psi_2=0$, $\psi_3(v_1,v_2)=\frac{1}{2}[v_1,v_2]$ and $\psi_4(v_1,v_2,v_3)=[v_1,v_3]+\frac{1}{6}\left[v_1,\left[v_1,v_2\right]\right]$. Therefore, if we write $v=\left(\begin{array}{c} v_k\\ \vdots\\ v_1\\ \end{array}\right)\in\mk{n}^+_1$, it holds: $$\begin{array}{rcl} \left(\begin{array}{c} k(g\cdot x)_k+\psi_k\left((g\cdot x)_1,\dots,(g\cdot x)_{k-1}\right)\\
\vdots\\ 3(g\cdot x)_3+\frac{1}{2}\left[(g\cdot x)_1,(g\cdot x)_2\right]\\ 2(g\cdot x)_2\\ (g\cdot x)_1\\
\end{array}\right) & = & d_g(x)_1^{-1}n_g(x) \\
\end{array}.$$
It follows that $[E,g\cdot x]=d_g(x)_1^{-1}n_g(x)+\psi(g\cdot x)$, where $\psi=\left(\begin{array}{c} \psi_k\\\vdots\\\psi_1\\\end{array}\right)$. More precisely, the map $\ad{E}$ is bijective on $\mk{n}^+_1$, so that we may obtain $(g\cdot x)_1=\text{pr}_{\g_1}\left(d_g(x)^{-1}_1n_g(x)\right)$, and \prtt $2\pt n\pt k$, $$(g\cdot x)_n=\frac{1}{n}\left(d_g(x)_1^{-1}n_g(x)\right)-\frac{1}{n}\psi_n\left((g\cdot x)_1,\dots,(g\cdot x)_{n-1}\right).$$
\end{demo}
Thus, if $\g$ is 3-graded, we recover the formula of \cite{BN1}, and in the general case of a $(2k+1)$-graded Lie algebra, we note that the projective elementary group acts in the chart $\mk{n}^+_1$ by birational maps.

\

As mentionned in the introduction, in this paper, we only consider purely algebraic issues. Using the differential calculus over a general topological ring defined in \cite{BGN}, it is possible, under suitable assumptions on the subalgebras $\mk{n}^+_1$ and $\mk{n}^-_1$ and on the general Bergman operators, to construct a structure of smooth manifold on $X^+$, resp. $X^-$, modelled on $\mk{n}^+_1$, resp. $\mk{n}^-_1$. All algebraically defined bundles are then smooth bundles, the canonical kernel is a smooth section and the tangent bundle defined in Definition \ref{defTiX} coincides with the tangent bundle as a smooth manifold. Such a construction is made in \cite{BN2} for a 3-graded Lie algebra and announced in \cite{moi1} for a general $(2k+1)$-graded Lie algebra.\\
Moreover, by definition $T_{o^+}X^+=\g/\mk{n}_0^-\cong\mk{n}^+_1$, and there is a filtration of $\mk{n}^+_1$ given by $$0\subset\g_1\subset\g_1\oplus\g_2\subset\dots\subset\mk{n}^+_1,$$ which is stable by the action of the group $P^-$, so that we obtain a filtration of $T_{o^+}X^+$. By the action of the projective elementary group, we also obtain a filtration of each $T_xX^+$, for $x\in X^+$. This distribution of filtrations is called a \emph{generalized contact structure} \cite{CKR}. This structure seems to be a key feature to understand the generalized flag geometry better. That is why we intend to come back to this structure in subsequent work.

\bibliographystyle{halpha}
\bibliography{biblio2}

\begin{thebibliography}{CDMKR05}

\bibitem[AF99]{AllisFaul}
B.~N. Allison and J.~R. Faulkner.
\newblock Elementary groups and invertibility for {K}antor pairs.
\newblock {\em Comm. Algebra}, 27(2):519--556, 1999.

\bibitem[Ber02a]{B02q}
W.~Bertram.
\newblock Complex and quaternionic structures on symmetric spaces -
  {C}orrespondance with {F}reudenthal-{K}antor triple systems.
\newblock {\em Theory of {L}ie Groups and Manifolds, Sophia Kokyuroku in
  Math.}, 45:57--76, 2002.

\bibitem[Ber02b]{B02}
W.~Bertram.
\newblock Generalized projective geometries: general theory and equivalence
  with {J}ordan structures.
\newblock {\em Adv. Geom.}, 2(4):329--369, 2002.

\bibitem[BGN04]{BGN}
W.~Bertram, H.~Gl{\"o}ckner, and K.-H. Neeb.
\newblock Differential calculus over general base fields and rings.
\newblock {\em Expo. Math.}, 22(3):213--282, 2004.

\bibitem[BK09]{BK1}
W.~Bertram and M.~Kinyon.
\newblock {A}ssociative {G}eometries. {I}: {G}rouds, linear relations and
  {G}rassmannians.
\newblock {\em Prépublication 2009/13 de l'Institut Elie Cartan}, 2009, arXiv:
  math.RA/0903.5441.

\bibitem[BL08]{BL08}
W.~Bertram and H.~L{\"o}we.
\newblock Inner ideals and intrinsic subspaces of linear pair geometries.
\newblock {\em Adv. Geom.}, 8(1):53--85, 2008.

\bibitem[BN04]{BN1}
W.~Bertram and K-H. Neeb.
\newblock Projective completions of {J}ordan pairs. {I}. {T}he generalized
  projective geometry of a {L}ie algebra.
\newblock {\em J. Algebra}, 277(2):474--519, 2004.

\bibitem[BN05]{BN2}
W.~Bertram and K-H. Neeb.
\newblock Projective completions of {J}ordan pairs. {II}. {M}anifold structures
  and symmetric spaces.
\newblock {\em Geom. Dedicata}, 112:73--113, 2005.

\bibitem[Bul07]{Bul07}
M.~Buliga.
\newblock Dilatation structures. {I}. {F}undamentals.
\newblock {\em J. Gen. Lie Theory Appl.}, 1(2):65--95 (electronic), 2007.

\bibitem[CDMKR05]{CKR}
M.~Cowling, F.~De~Mari, A.~Kor{\'a}nyi, and H.~M. Reimann.
\newblock Contact and conformal maps in parabolic geometry. {I}.
\newblock {\em Geom. Dedicata}, 111:65--86, 2005.

\bibitem[Che09]{moi1}
J.~Chenal.
\newblock Generalized flag geometries and manifolds associated to short
  {$\mathbb Z$}-graded {L}ie algebras in arbitrary dimension.
\newblock {\em C. R. Math. Acad. Sci. Paris}, 347(1-2):21--25, 2009.

\bibitem[Che10]{moi2}
J.~Chenal.
\newblock Structures géométriques liées aux algèbres de lie graduées.
\newblock {\em PhD Thesis}, 2010.

\bibitem[Kam89]{Kam}
Noriaki Kamiya.
\newblock A structure theory of {F}reudenthal-{K}antor triple systems. {III}.
\newblock {\em Mem. Fac. Sci. Shimane Univ.}, 23:33--51, 1989.

\bibitem[Kan70]{Kan70}
I.~L. Kantor.
\newblock Graded {L}ie algebras.
\newblock {\em Trudy Sem. Vektor. Tenzor. Anal.}, 15:227--266, 1970.

\bibitem[Kan72]{Kan72}
I.~L. Kantor.
\newblock Some generalizations of {J}ordan algebras.
\newblock {\em Trudy Sem. Vektor. Tenzor. Anal.}, 16:407--499, 1972.

\bibitem[Loo75]{Loos75}
O.~Loos.
\newblock {\em Jordan pairs}.
\newblock Lecture Notes in Mathematics, Vol. 460. Springer-Verlag, Berlin,
  1975.

\bibitem[Mey70]{Mey}
K.~Meyberg.
\newblock Jordan-{T}ripelsysteme und die {K}oecher-{K}onstruktion von
  {L}ie-{A}lgebren.
\newblock {\em Math. Z.}, 115:58--78, 1970.

\bibitem[Rub94]{Rub}
Hubert Rubenthaler.
\newblock Les paires duales dans les alg\`ebres de {L}ie r\'eductives.
\newblock {\em Ast\'erisque}, (219):121, 1994.

\bibitem[Upm85]{Up}
H.~Upmeier.
\newblock {\em Symmetric {B}anach manifolds and {J}ordan {$C\sp
  \ast$}-algebras}, volume 104 of {\em North-Holland Mathematics Studies}.
\newblock North-Holland Publishing Co., Amsterdam, 1985.
\newblock Notas de Matem{\'a}tica [Mathematical Notes], 96.

\end{thebibliography}

\end{document}